\documentclass[11pt]{article}
\usepackage[english]{babel}
\usepackage[cp850]{inputenc}
\usepackage[dvips]{graphicx}
\usepackage{amsmath,amsfonts,amsthm,amssymb}
\usepackage[usenames, dvipsnames]{color}
\usepackage{fancyhdr}
\usepackage{stmaryrd}
\usepackage[colorlinks=true,citecolor=red,linkcolor=blue,urlcolor=blue,pdfstartview=FitH]{hyperref}
\usepackage{dsfont}
\usepackage{xcolor}
\usepackage{epsfig}

\font\tenmath=msbm10 scaled 1200

\font\sevenmath=msbm7 scaled 1200
\font\fivemath=msbm5 scaled 1200

\newfam\mathfam \textfont\mathfam=\tenmath
\scriptfont\mathfam=\sevenmath \scriptscriptfont\mathfam=\fivemath

\def\R{{\mathbb R}}
\def\N{{\mathbb N}}
\def\E{{\mathbb E}}

\def\P{{\mathbb P}}

\newtheorem{Thm}{Theorem}[section]
\newtheorem{Lem}{Lemma}[section]
\newtheorem{Pro}{Proposition}[section]

\newtheorem{Dfn}{Definition}[section]

\newfam\mathfam \textfont\mathfam=\tenmath
\scriptfont\mathfam=\sevenmath \scriptscriptfont\mathfam=\fivemath

\def \^#1{\if#1i{\accent"5E\i}\else{\accent"5E#1}\fi}

\title{\bf  Pointwise 
Convergence of the Lloyd algorithm in higher dimension}

\author{
\textsc{Gilles Pag\`es} \thanks{UPMC,  Laboratoire de Probabilit\'es et Mod\`eles al\'eatoires, UMR~7599, case 188, 4, pl. Jussieu, F-75252 Paris Cedex 5, France. E-mail: \texttt{gilles.pages@upmc.fr}}  
\quad \textsc{Jun YU}~\thanks{UPMC, Laboratoire de Probabilit\'es et Mod\`eles al\'eatoires, E-mail: \texttt{jun.yu@upmc.fr}}
}

\begin{document}
\maketitle
\begin{abstract} We establish the pointwise convergence of the iterative Lloyd algorithm, also known as $k$-means algorithm,  when the quadratic quantization error of  the starting grid (with size $N\ge 2$) is lower than the minimal quantization error with respect to the input distribution is lower at level $N-1$. Such a protocol is known as the  splitting method and allows for convergence even when the input distribution has an unbounded support. We also show under very light assumption that the resulting limiting grid still has full size $N$. These results are obtained without continuity  assumption on  the input distribution. A variant of the procedure taking advantage of the asymptotic of the optimal quantizer radius is proposed which always  guarantees the boundedness of the iterated grids. 
\end{abstract}


\section{Introduction} 

A Centroidal Voronoi Tessellation (CVT) with respect to a probability (or mass) distribution is a   Voronoi
tessellation of a  set of (generating) points in $\R^d$ (centers of mass)  such that each  generating point is the centroid of
its  corresponding Voronoi region with respect to this density function. This definition can be extended to more general probability measures, typically those assigning no mass to hyperplanes to avoid ambiguity on the boundaries of the Voronoi regions. CVTs   enjoy very natural optimization
properties, especially in connection with vector quantization (see further on)  which makes them very popular in various scientific and engineering applications including art design, astronomy, clustering, geometric modeling, image and
data analysis, resource optimization, quadrature design, sensor networks, and numerical solution of partial differential equations.

\medskip For modern applications of the CVT concept in large-scale scientific and engineering problems, it is important to develop robust and efficient algorithms for constructing CVTs in various settings. Historically, a number of algorithms have been studied and widely used. However, the pioneering  contribution is undoubtedly the procedure first developed in the 1960s at Bell Laboratories by S. Lloyd. It remains so far, in its randomized form,  one of the most popular methods due to its effectiveness and simplicity.

\medskip Let us begin  with a more detailed description of the CVT. First assume that the probability distribution, say $\mu$, on $(\R^d, {\cal B}or(\R^d))$, has a support included in a closed convex set with non empty interior denoted $\mathbf{U}$ of $\R^d$. Also, note that, up to a reduction of the dimension $d$, one may always assume that $\mathbf{U}$ has a nonempty interior. 

A Voronoi diagram (or partition)  of $\mathbf{U}$ refers to a Borel partition  $(C_i(\Gamma))_{1\le i\le N}$ of $\mathbf{U} \subset \mathbb{R}^d$   induced by a set $\Gamma= \{x_i,\, 1\le i\le N\} \subset \mathbf{U}$ of $N$ given {\em generating points} or {\em  Generators}   (the notation $\Gamma$ also refers to the application to numerics where the set of generators is also called a {\em grid}). For every  $i\!\in \{1,\ldots,N\}$, the {\em Voronoi region} (or {\em cell}) $C_i(\Gamma)$ satisfies
$$
C_i(\Gamma)\subset \left\{\xi \!\in \mathbf{U}\,:\, |\xi-x_i|\le \min _{1\le j\le N} |\xi-x_j|\right\}
$$ 
where $|\,.\,|$ denotes the canonical Euclidean norm on $\R^d$. Then
\[
 \left\{\xi \!\in \mathbf{U}\,:\, |\xi-x_i|< \min _{1\le j\le N} |\xi-x_j|\right\}= \stackrel{\circ}{C}_i(\Gamma)\subset C_i (\Gamma)\subset \overline{C}_i (\Gamma)=  \left\{\xi \!\in \mathbf{U}\,:\, |\xi-x_i|\le \min _{1\le j\le N} |\xi-x_j|\right\}
\] 
so that the  $C_i(\Gamma)$ have convex interiors and closures. The family of closures is also known as {\em Voronoi tessellation} of $\mathbf{U}$ induced by $\Gamma$ and the $\overline{C}_i(\Gamma)$, $i=1,\ldots,N$ are called {\em tessels}). Furthermore they have a polyhedral structure, in particular their boundaries are contained in $\displaystyle \cup_{i \neq i}H_{ij}$ where $H_{ij}\equiv \frac{x_i+x_j}{2} +  \Big(\frac{x_i-x_j}{|x_i-x_j|}\Big)^{\perp}$ is the median hyperplane of $x_i$ and $x_j$.   Of course,  a notion of Voronoi regions can be defined with respect to  any norm $N$ on $\R^d$ but the above (polyhedral) convexity  properties fail  (see $e.g.$~\cite{Fou}, chapter~1) for non Euclidean norms.

We will often assume that  $\mu$ is {\em strongly continuous} in the sense that it assigns no mass to hyperplanes (so is the case if $\mu$ is absolutely continuous $i.e.$ $\mu(d\xi) =\rho(\xi)d\xi$ where $\rho$ is a probability density function defined on $\R^d$ whose support is contained in $\mathbf{U}$). Then  the boundaries of the Voronoi regions  are $\mu$-negligible so that we can define in a unique way the {\em centroids} $x^*_i$, $i=1,\ldots,N$ of the Voronoi regions by setting
 \begin{equation}\label{eq:x*} 
 x_i^* =\left\{\begin{array}{ll}
\displaystyle   \frac{\int_{C_i} \xi \mu(d\xi)}{\mu(C_i)} & \mbox{if } \mu(C_i)>0,\\
  x_i & \mbox{if } \mu(C_i)=0,
  \end{array}\right. \; i=1,\ldots,N.
 \end{equation}
 Note that, owing to the convexity of the Voronoi cells $C_i$ and the finiteness of the measure $\mu$, one has    $x^*_i\!\in \overline C_i$ (closure in $\mathbf{U}$) for every $i\!\in \{1,\ldots,N\}$. From a more probabilistic point of view, if  $X$ denotes an $\R^d$-valued  random vector with distribution $\P_{_X}=\mu$, then  (with an obvious convention when $\P(X\!\in C_i)=0$) 
 \[
 x^*_i = \E \big(X\,|\, X\in C_i\big), \; i=1,\ldots,N.
 \]
This naturally leads to the definition of a CVT which is but a  Voronoi tessellation  whose generators $x_i$ are   the centroids of their respective Voronoi regions. With the notation given above, the Lloyd algorithm for constructing CVTs can be described more precisely by the following procedure.

\medskip The paradigm of Lloyd's algorithm is to consider the definition of CVT as a fixed point equality for the so-called {\em Lloyd map} $T^{\mu}_{_N}$ defined on  the set of $\mathbf{U}$-valued  grids $\Gamma$  with at most $N$ values by~\eqref{eq:x*}, $i.e.$ 
\[
T^{\mu}_{_N}(\Gamma) = \{x^*_i ,\; i=1,\ldots,N\}\quad\mbox{if} \quad \Gamma= \{x_i ,\; i=1,\ldots,N\}\subset \mathbf{U}.
\]
As mentioned above,  $T(\Gamma)_i\!\in \overline C_i(\Gamma)$ since the Voronoi tessels are convex and $\mu$ is probability distribution. 

Note that, furthermore,  if ${\rm supp}(\mu)= \mathbf{U}$ or if $\mu$ is contain,ious (assigns no mass to hyperplanes) then, a supporting hyperplane argument shows that  $T(\Gamma)_i\in \,  \stackrel{\circ}{C}_i(\Gamma)$ (interior in $\mathbf{U}$) for every $i=1,\ldots,N$  (see further on~Lemma~\ref{lem:nomerge}, see also~\cite{Fou}, p.22). In particular, $T(\Gamma)$ and $\Gamma$ have the same size  $N$.

\medskip  The Lloyd algorithm is simply the formal fixed point search procedure for the Lloyd map $T_{_N}$ starting from a given grid $\Gamma^{(0}= \{x^{(0)}_i,\, i=1,\ldots,N\}$ of full size $N$ $i.e.$
\[
\Gamma^{(k+1)}= T^{\mu}_{_N} (\Gamma^{(k)}), \; k\ge 0.
\]

\noindent {\bf Algorithm 1 (Lloyd's  algorithm for computing CVTs)}:

\smallskip
\noindent $\rhd$ {\em Inputs}:
\begin{itemize}
\item  $\mathbf{U}$, the domain of interest; 
\item  $\mu$ a probability distribution supported by  $\mathbf{U}$;
\item  $\Gamma^{(0)} = \{x_i^{(0)} ,\,i=1,\ldots,N \}\subset \mathbf{U}$, the initial set of $N$ generators.
\end{itemize}

\noindent $\rhd$ {\em Pseudo-script}:

\smallskip
Formally, at the $k$\textsuperscript{th} iteration, one has to proceed as follows:
\begin{enumerate}
\item Construct the Voronoi tessellation $\{C_i(\Gamma^{(k)}),\; i=1,\ldots,N\}$ of $\mathbf{U}$ with the grid of generators $\Gamma^{(k)} = \{x_i^{(k)} ,\; i=1,\ldots,N\}$.
\item Compute the $\mu$-centroids of $\{C_i(\Gamma^{(k)}),\; i=1,\ldots,N\}$ as the new grid  of generators $\Gamma^{(k+1)} = \{x_i^{(k+1)},\; i=1,\ldots,N\}$.
\end{enumerate}

\noindent  $\rhd$  Repeat the iteration above until some stopping criterion is met to provide a grid of generators as close as possible of a $\mu$-centroid. 
%

\smallskip
\noindent $\rhd$  end.

\medskip
In $1$-dimension, Kieffer has proved in~\cite{KIE} that $T_{_N}$ is contracting if $\mu$ has a $\log$-concave density over a compact interval so that  only one $\mu$-centroid with $N$ points exists for such distribution and the above procedure converges exponentially fast toward it.  See also,  more recently a convergence result in~\cite{Du2}.

\medskip
In practice, these two steps become intractable in higher dimension  by analytic or even deterministic approximation methods, say when $d\ge 3$ or $4$ (see however the website {\em QHull}: {www.qhull.org}). So this theoretical  procedure has to be replaced for numerical purpose by a randomized  version in which:

\smallskip
-- Step 1 is replaced by a systematic nearest neighbour search of simulated random $\mu$-distributed vectors.

\smallskip 
-- Step 2 is replaced by a {\em Monte Carlo estimation} of both terms of the ratio which  define the Lloyd map. 

\medskip
In the community of data analysis, note that when $\mu=\frac 1M \sum_{m=1}^M\delta_{\xi_m}$ is the empirical measure of a  $\mathbf{U}$-valued data set $(\xi_m)_{m=1,\ldots,M}$, it is still possible to define and compute the Lloyd map (using appropriate conventions like $e.g.$ random allocation of points lying on the boundary of (closed) Voronoi tessels). In such a case the Lloyd procedure is known as the Forgy algorithm or the batch-$k$-means procedure. When the data set is so huge that a   uniform sampling (of size $M$) of the dataset is necessary at each iteration, the procedure is known as the {$k$-means} procedure.

\smallskip In this paper we will focus on the  converging properties  of the theoretical (or {\em batch} in the data-mining community) Lloyd procedure,  prior to any randomization or approximation, although we are aware that in higher dimension for continuous distributions $\mu$, it is a {\em pseudo-algorithm}. So far we have presented the Lloyd procedure in an intrinsic manner. In fact Lloyd's algorithm is deeply connected with the theory of Optimal Vector Quantization of probability distribution. This connection turns out to provide very powerful tools to investigate its convergence properties. It is also a major field of application  when trying to compute with a sharp accuracy optimal quantizers of  simulatable distribution arising  in the design of numerical schemes for solving nonlinear problems (optimal stopping problems, (possibly Reflected) Backward Stochastic Differential Equations, Stochastic Control, etc, see~\cite{BAPA, PAPHPR}).

\medskip 
Quantization is a way to discretize the path space of a random phenomenon: a random vector in finite dimension (but also stochastic process in infinite dimension viewed as a random variable taking values in its path space which we will not investigate in this paper). We consider here a random vector $X$ defined on a probability space $(\Omega, \mathcal{A}, \mathbb{P})$ taking its values in $\R^d$ equipped with its Borel $\sigma$-field ${\cal B}or(\R^d)$.

\medskip It is convenient for what follows to introduce a few notions and results about vector quantization and its (mean quadratic) optimization. It makes a connection between CVTs and stochastic optimization, gives a rigorous meaning to the notion of ``goodness" of a CVT. Optimal vector quantization goes back to the early 1950's in the Bell laboratories and have been developed for the optimization of signal transmission.
 
Let $X:(\Omega, \mathcal{A}, \mathbb{P}) \rightarrow \R^d$  be a square integrable random vector  ($i.e.$ $\E|X|^2<+\infty$) or equivalently $X\! \in L_{\R^d}^2(\mathbb{P})$. Assume that its distribution $\mu = \P_{_X}$ is included in~$\mathbf{U}$ (defined as above). The terminology  {\em $N$-quantizer}  (or a {\em quantizer at level} $N$) is assigned to  any $\mathbf{U}$-valued subset  with cardinality $N$. When used in a numerical framework, it is also known as {\em quantization grid}. 
\[ 
\Gamma := \{x_1, x_2, \cdots , x_N \} \subset \mathbf{U}.
\]
The cardinality of $\Gamma$ is $N$. In numerical applications, $\Gamma$ is also called a (quantization) {\em grid}. It is the set of genrators of its (borel) Voroni regions $(C_i(\Gamma))_{1\le i\le N}$. Then can discretize $X$ in pointwise way by $q(X)$ where  $q$: $\R^d \rightarrow \Gamma $  is a Borel function. Then we get 
\[ 
\forall \omega\! \in\Omega, \; |X(\omega) - q(X(\omega))| \ge {\rm dist}(X(\omega), \Gamma) = \min_{1\le i \le N} |X(\omega) - x_i| 
\]
so that the best pointwise approximation of $X$ is provided by considering any (Borel) {\em nearest neighbour projection} $q={\rm Proj}_{\Gamma}$   associated with the  Voronoi tessellation $(C_i(\Gamma))_{1\le i\le N}$ by  setting
\[ 
\operatorname{Proj}_{\Gamma} (\xi) = \sum_{i=1}^N x_i \mathbf{1}_{C_i(\Gamma)}(\xi) ,\; \xi\!\in \R^d.
\]
It is clear that such a  projection is in one-to-one correspondence with the Voronoi partitions (or diagrams) of $\R^d$ induced by $\Gamma$.
%
These projections only differ on the boundaries of the Voronoi cells $C_i(\Gamma)$ so that, as soon as $\mu = \P_{_X}$  is strongly contoinuous, these neratest neighbour projections are all $\P_{_X}$ -$a.s.$ equal.   We define a {\em Voronoi $N$-quantization} of $X$ (or at level $N$)  by setting for every $\omega \!\in \Omega$,
\[ 
\widehat{X}^{\Gamma}(\omega) := \operatorname{Proj}_{\Gamma} (X(\omega)) = \sum_{i=1}^N x_i \mathbf{1}_{C_i(\Gamma)}(X(\omega)).
\]
Thus for all $\omega \in \Omega$,
\begin{equation}\label{eq:VoroGeo} 
|X(\omega) - \widehat{X}^{\Gamma}(\omega) | ={\rm dist}(X(\omega), \Gamma) = \min_{1\le i \le N} |X(\omega) - x_i|. 
\end{equation}
We will call $\widehat{X}^{\Gamma}$ a Voronoi $\Gamma$-quantization of $X$ or, in short, a quantization of $X$. 

\medskip The mean quadratic quantization error is then defined by
\[ 
e(\Gamma, X) = \|X - \widehat{X}^{\Gamma} \|_2 = \sqrt{\mathbb{E}\left( \min_{1\le i \le N} |X - x_i|^2 \right)}
\]
where $\|\cdot\|_2$ is the norm in $L_{\R^d}^2(\mathbb{P})$. The distribution of $\widehat{X}^{\Gamma}$ as a random vector is given by the $N$-tuple $\left(\mathbb{P}(X \in C_i(\Gamma) ) \right)_{1\le i \le N}$. This distribution clearly depends on the choice of the Voronoi partition.

\medskip We naturally wonder whether it is possible to design some optimally fitted grids to a given distribution $\mu = \mathbb{P}_X$ \textit{i.e.} which induces the lowest possible mean quadratic  quantization error among all
grids of size at most $N$. This optimization problem, known as {\em the optimal quantization problem at level $N$}, reads as follows:
\[ 
e_N(X) := \inf_{\Gamma \subset \mathbb{R}^d, \operatorname{Card}(\Gamma)\le N} e(\Gamma, X ) 
\]
By introducing the energy function or  {\em  distortion value function}
\[ 
\mathcal{G} : \, \big(\mathbb{R}^d\big)^N \longrightarrow  \mathbb{R}_+
\] 
\[ 
x=(x_1, x_2, \cdots , x_N) \longmapsto \mathbb{E}\left( \min_{1\le i \le N} |X - x_i|^2 \right)
\]
the optimization problem also reads 
\[
e_N(X)  = \inf_{x \in (\mathbb{R}^d)^N} \sqrt{\mathcal{G} (x)}
\]
since the value of $\mathcal{G}$ at an $N$-tuple $x=( x_1, x_2, \cdots, x_N)\!\in (\R^d)^N$ only depends on its value grid $\Gamma=  \Gamma_x=\{x_1,\ldots,x_N\}$ of size at most $N$ of the $N$-tuple (in particular $\mathcal{G}$ is a symmetric function). We will make occasionally the abuse of notation consisting in denoting $\mathcal{G}(\Gamma)$ instead of $\mathcal{G}(x)$.
%
%

One proves  (see $e.g.$~\cite{Cuesta-albertos, Fou, PAG0}) that there always exists at least one  optimal $N$-point grid $\Gamma_{_N}^* = \{x_1^*, x_2^*, \cdots , x_N^*\} \subset \mathbb{R}^d $ with cardinal $N$  such that $e_N(X)=  \sqrt{\mathcal{G} (\Gamma_{_N}^*)} $. If the support of $\P_{_X}$ has at least  $N+1$ elements ($e.g.$ because it is infinite), then $\Gamma_{_N}^*$ has full size $N$. Furthermore, $\Gamma^*_{_N}\subset \mathbf{U}$; this last claim strongly relies on  the Euclidean feature of  the norm on $\R^d$: if $\Gamma^*_{_N} \,/ \hskip -0.4 cm  \subset \mathbf{U}$, then the projection of the elements of $\Gamma^*_{_N}$ on  the closed (nonempty) convex $\mathbf{U}$  strictly reduces  the mean quadratic quantization error   (see $e.g.$~\cite{PAG0, Fou, GRLUPA3}). 
Note that this existence result does not require $\mu= \P_{X}$ to be strongly continuous. In fact, even if $\mu$ has an atomic component, it is shown in~\cite{Fou} (see Theorem~4.2, p.38) that $\mu\big(\cup_i \partial C_i(\Gamma^*_{_N})\big)=0$.

\smallskip
 Furthermore, the function $\mathcal{G}$ is differentiable on $(\R^d)^N$ at every $N$-tuple $x=(x_1,\ldots, x_{_N})$ such that $x_i\neq x_j$, $i\neq j$ and $\displaystyle \P\big( X \!\in \cup_{1\le i\le N} \partial C_i(x)\big) =0$ and its gradient is given by 
\begin{equation}\label{eq:gradG}
\nabla \mathcal{G}(x)= \frac 12 \E\left(\mbox{\bf 1}_{\{X\!\in C_i(x)\}}( x_i-X)\right)
\end{equation}
where $(C_i(x))_{1\le i\le N}$  denotes any Voronoi diagram of $\{x_1,\ldots,x_{_N}\}$ (or $x$ with, once again, an obvious abuse of notation).
In particular, if $x^*=(x^*_1,\ldots,x^*_{_N})$ is  an optimal quadratic quantizer, one shows (see~\cite{Fou}, Theorem~4.2 p.38) that $\P\big( X \!\in \cup_i \partial C_i(x^*)\big) =0$ even if $\P_{_X}$ assigns mass to (at most countably many)  hyperplanes. As a consequence, 
\[
\nabla \mathcal{G}(x^*)= 0\quad\mbox{ or equivalently} \quad T^{\P_{_X}}_N (x^*)=x^*
\]
where $T^{\P_{_X}}_N$ stands for the Lloyd map related to the distribution $\mu = \P_{_X}$ of the random vector $X$. The equivalence is a straightforward consequence of~\eqref{eq:gradG}). Hence any optimal quantizer induces a CVT for the distribution of $X$. In reference to the fact that such an $N$-quantizer is a zero of a gradient, this property is also known as {\em stationarity} for the $N$-tuple $x^*$ itself.

Unfortunately, the converse is not true since $\mathcal{G}$ may have many local minima, various types of saddle points (and a ``pin"  behaviour on affine manifolds induced by clusters of stuck components). This phenomenon becomes more and more intense as $d$ grows. However, it makes a strong connection between search for optimal quantizers and Lloyd's algorithm as described above. And there is no doubt that what practitioners are interested in are the optimal quantizers rather than any ``saddle" stationary quantizers. One can also derive a stochastic gradient algorithm from the representation of $\nabla \mathcal{G}$ as an expectation of a computable function of the quantizer $x$ and the random vector $X$. This second approach leads to a   stochastic optimization procedure, a stochastic gradient descent to be more precise,  known as the {\em Competitive learning vector Quantization} algorithm ($CLVQ$) which has also been extensively investigated (see among others~\cite{PAG0}).

\medskip The paper is organized as follows: in Section~\ref{sec:2} we establish the convergence of the Lloyd procedure at level $N$ under some natural assumptions on the probability distribution (at least for numerical probability purpose) but assuming that the starting quantizer (or generators) induces a lower quantization error than the lowest quantization error at level $N-1$.  In Section~\ref{sec:3}, we propose a modified Lloyd's procedure, inspired by recent results on the asymptotics of the ``radius" of optimal quantizers at level $N$ as $N\to +\infty$, to overcome partially this constraint on the starting grid.

\medskip
\noindent {\sc Notations:} $\bullet$  ${\rm supp}(\mu)$ denotes the support of the Borel probability measure $\mu$ on $\R^d$. $\lambda_d$ denotes the Lebesgue measure on $\R^d$. 

\noindent  $\bullet$ $(.\,,.)$ denotes the canonical inner product  on $\R^d$. $B(x, r)$, $r\!\in \R_+$, denotes the canonical Euclidean ball centered at $x\!\in \R^d$ with radius $r>0$. 

\noindent  $\bullet$ $|A|$ denotes the cardinality of set $A$.

\section{Convergence analysis of Lloyd's algorithm with unbounded inputs}\label{sec:2} 
\subsection{The main result}\label{subset:2.1}
Owing to both simplicity   and efficiency of   of practical implementations of Lloyd's algorithm in various fields of applications, it is important to study its convergence as it has been carried out, at least partially,  in~\cite{PAG0} for its ``counterpart" in the world of Stochastic Approximation, the recursive stochastic gradient descent attached to the above gradient of the distortion function $\mathcal{G}$. This procedure is also known as $CLVQ$ (for Competitive Learning Vector Quantization algorithm). In fact, as concerns the convergence properties of Lloyd's algorithm, many investigations have already been carried out . Thus, as mentioned in the introduction, true convergence for $\log$-concave densities has been established in~\cite{KIE} whereas  global ``weak"  convergence has been proved   in a one dimensional setting (see \cite{Du2}). However, there are not many general mathematical results on the convergence analysis for  
distributions on multi-dimensional spaces, especially when {\em the support of the distribution $\mu$ of interest  is not bounded}. 

It is convenient to rewrite the iterations of the Lloyd algorithm in a more probabilistic form, using quantization formalism (with a generic notation for the grids: $\Gamma^{(k)}= \big\{x^{(k)}_1,\ldots,x^{(k)}_{_N}\big\}$, $k\ge 0$). Let $\Gamma^{(0)} \subset \R^d$. For every $k\ge 0$,  
\[
 \left\{\begin{array}{lll} (I) \hbox{\em Centroid updating}: &\hskip -1 cm \widetilde X^{(k+1)}= \E\big( X\,|\, \widehat X^{\Gamma^{(k)}}\big),\; \Gamma^{(k+1)}&= \widetilde X^{(k+1)}(\Omega)\\
 &&  = \left\{\frac{\E\big(\mbox{\bf 1}_{\{\widehat X^{\Gamma^{(k)}}= x^{(k)}_i\}}X\big)}{\P(\widehat X^{\Gamma^{(k)}}= x^{(k)}_i)},\,i=1,\ldots,N\right\}\\
(II) \hbox{\em Voronoi  cells re-allocation}:&  \widehat X^{\Gamma^{(k+1)}} \leftarrow \widetilde X^{(k+1)}.&
 \end{array}\right. 
\] 
with the following {\em degeneracy convention}:

\medskip
\centerline{ if $\P(\widehat X^{\Gamma^{(k)}}= x^{(k)}_i) = \P( X\!\in C_i(\Gamma^{(k)}) )=0$ then $\widehat X^{\Gamma^{(k+1)}}=x^{(k)}_i$.}
 

\medskip We need now to formalize in a more precise way what can be a {\em consistent connection} between the sequence of iterated grids $(\Gamma^{(k)})_{k\ge 0}$ in the Lloyd procedure  and the $N$-tuples  that can be associated to their values.
\begin{Dfn} Let $(\Gamma^{(k)})_{k \ge 0}$ be a sequence of iterates of the Lloyd procedure where $|\Gamma^{(0)}|= N$.  A sequence of $N$-tuples 
$(x^{(k)})_{k\le  0}$  is a {\em  consistent representation} of the sequence $\Gamma^{(k)}$ if
\begin{itemize}
\item[(i)] $\Gamma^{(k)}= \{x^{(k)}_i,\, i=1,\ldots,N\}$,
\item[(ii)] For every integer $k\ge 0$ and every $i\!\in \{1,\ldots,N\}$, $x^{(k+1)}_i $ is the centroid of the cell of $x^{(k)}_i $ $i.e.$
\[
x^{(k+1)}_i = \E\big(X\,|\, X\!\in C_i(\Gamma^{(k)})\big)= \frac{\E\big(\mbox{\bf 1}_{\{\widehat X^{\Gamma^{(k)}}= x^{(k)}_i\}}X\big)}{\P(\widehat X^{\Gamma^{(k)}}= x^{(k)}_i)}
\]
\end{itemize}
(still with the above degeneracy convention). 
\end{Dfn}

There are clearly $N!$ consistent representations of a sequence of Lloyd iterates, corresponding to the possible numbering of $\Gamma^{(0)}$. But then this numbering is frozen as $k$ increases. It is also clear that roundedness, convergence (in $(\R^d)^N$) of such consistent representations does not depend on the selected representation. So is true for a possible limit $x^{(\infty)}$ of such sequences since $\Gamma^{(\infty)}$ will not depend on the selected representation. However on may have $|\Gamma^{(\infty)}|<N$ in case of an asymptotic  merging of some of the components of the representation. One checks that under the assumptions  we make ($\mu$ continuous, or convex support or splitting assumption on $\Gamma^{(0)}$) no merging   occurs at finite range. 

\medskip Throughout the paper $(x^{(k)})_{k\ge 0}$ will always denote a consistent representation of the sequence of iterates $\Gamma^{(k)})_{k\ge 0}$. These remark lead naturally to the following definition

\begin{Dfn}[Convergence of iterated grids] \label{ConvGrid}We will say that $\Gamma^{(k)}\to \Gamma^{(\infty)}$ (converges) as $k\to+\infty$ if there exists a constant representation $(x^{(k)})_{k\ge 0}$ converging in $(\R^d)^N$ toward $x^{(\infty)}$ such that $\Gamma^{(\infty)}= \{x^{(k)}_i, \, i=1,\ldots,N\}$.
\end{Dfn}

Several specific results established in what follows are  known in the literature, but we chose to provide all  proofs for self-completeness of the paper and reader's convenience.

%

Let us first recall a basic fact which is at the origin of the efficiency of the Lloyd algorithm.

\begin{Lem} \label{Lem:seminal}The iteration of Lloyd's algorithm makes the quantization error decrease. 
\[
e(\Gamma^{(k+1)}, X ) \le  e(\Gamma^{(k)}, X )  .
\] 
Furthermore, $e(\Gamma^{(k+1)}, X ) <  e(\Gamma^{(k)}, X )  $ as long as $ \widehat X^{\Gamma^{(k)}}\neq  \mathbb{E}( X | \widehat{X}^{\Gamma^{(k)}} ) $ with positive $\P$-probability. Conversely, if $ \widehat X^{\Gamma^{(k_0)}} =   \mathbb{E}( X | \widehat{X}^{\Gamma^{(k_0)}} ) $ $\P$-$a.s.$ for an integer $k_0\!\in \N$, then $ \widehat X^{\Gamma^{(k)}}= \widehat X^{\Gamma^{(k_0)}}$ $\P$-$a.s.$ for every $k\ge k_0$.
\end{Lem}
\noindent\textbf{Proof}:
By its very definition, $\widetilde{X}^{(k+1)}$ is the best approximation of $X$ among the $\sigma(\widehat{X}^{\Gamma^{(k)}})$-measurable functions, including $\widehat{X}^{\Gamma^{(k)}}$ itself. Thus
\begin{equation}\label{eq:step1}
 \| X - \widetilde{X}^{(k+1)} \|_2 = \| X -  \mathbb{E}( X | \widehat{X}^{\Gamma^{(k)}} ) \|_2  \le \| X - \widehat{X}^{\Gamma^{(k)}} \|_2= e(\Gamma^{(k)}, X )
\end{equation}
with equality if and only if $ \widehat{X}^{\Gamma^{(k)}} = \mathbb{E}( X | \widehat{X}^{\Gamma^{(k)}} )$.   Note that, if so is the case, $\Gamma^{(k+1)} = \Gamma^{(k)}$ since $\widetilde X^{(k+1)}= \widehat X^{\Gamma^{(k)}}$. On the other hand, by~\eqref{eq:VoroGeo},  $\widehat X^{\Gamma^{(k+1)}}$ is the best  approximation of $X$ among all $\Gamma^{(k+1)}$-valued random vectors since
\[ 
e(\Gamma^{(k+1)}, X )=\| X - \widehat{X}^{\Gamma^{(k+1)}} \|_2 = \| {\rm dist}(X, \Gamma^{(k+1)})\|_2  \le \| X - \widetilde{X}^{(k+1)} \|_2
\] 
and if $\Gamma^{(k+1)} = \Gamma^{(k)}$, this inequality holds as an equality.
$\qquad\Box$

%
%
%

\bigskip This seminal property of Lloyd's algorithm is striking not only by its simplicity. It is also at the origin of its success. Morally speaking, it suggests a convergence toward a stationary --~and hopefully optimal or at least locally optimal~-- quantizer of the distribution $\mu$ of $X$. In fact, things are less straightforward, at least from a theoretical point of view since this property provides absolutely no information on the boundedness of the sequence of grids $(\Gamma^{(k)})_{k\ge 0}$ generated by the procedure, although it is a crucial property the way toward convergence.

\medskip The main result of this paper is the following theorem.

 \begin{Thm} \label{thm:main} Let $X:(\Omega,{\cal A}, \P)\to \R^d$ be a square integrable random vector with a distribution $\mu$ having a  convex support $\mathbf{U}$   assigning no mass to hyperplanes.  Let $\Gamma^{(0)}\subset \mathbf{U}$ with size $|\Gamma^{(0)}|=N$. Then, all the iterates $\Gamma^{(k)}=\{x^{(k)}_i,\, i=1,\ldots,N\}$ have full size $N$ and are $\mathbf{U}$-valued.

\smallskip
 \noindent $(a)$ If the sequence $\big(\Gamma^{(k)}\big)_{k\ge 1}$ is bounded (or equivalently kits consistent representations  $(x^{(k)})_{k\ge 0}$), then
 \[
 \liminf_k \min_{1\le i,j\le N}|x^{(k)}_i-x^{(k)}_j|>0
 \] 
 and there exists $\ell_\infty\! \in \big[0, e_{N}(\Gamma^{(0)}(X))\big)$ and a connected component  ${\cal C}_{\ell_\infty, \mathbf{U}}$ of  $\displaystyle \Lambda_{\infty} :=\big\{x\!\in\mathbf{U}^N, \, e_N\big(\{x_i,i=1,\ldots,N\},X\big)=\ell_{\infty}, \; \mathcal{G}(x)=0\big\}$ such that 
\[
{\rm dist}\big(x_i^{(k)}, {\cal C}_{\ell_\infty, \mathbf{U}}\big)\to 0\quad \mbox{as}\quad k\to +\infty.
\]
In particular, if $\Lambda_{\infty}$ is locally finite ($i.e$ is reduced to finitely many points on each compact set), then $\Gamma^{(k)}\to \Gamma^{(\infty)}\!\in \Lambda_{\infty}$.  

\smallskip
\noindent $(b)$ {\em Splitting method:} If furthermore $e(\Gamma^{(0)},X) \!  \in (e_{N}(X), e_{N-1}(X)]$, then the sequence $\big(\Gamma^{(k)}\big)_{k\ge 1}$   is always $\P$-$a.s.$ bounded and $\ell_{\infty}\!\in  [e_{N}(X), e_{N-1}(X))$.
 \end{Thm}

 \noindent {\bf Remarks.} $\bullet$ The result in claim $(a)$ does not depend on the original numbering of $\Gamma^{(0)}$ $i.e.$ on the selected order (among $N!$) selected to define $x^{(0)}$ and the then the sequence $(x^{(k)})$. In particular in case of true convergence of the sequence $x^{(k)}$, all its permutations do converge as well toward the corresponding permutations of $x^{(\infty)}$ so that one can, by an abuse of notion write that $\Gamma^{(k)}\to \Gamma^{(\infty)}$. A direct approach based on a formal notion of set  convergence  is also possible but would be of no help in practice.  
 
 \smallskip 
 \noindent $\bullet$ Nothing ensures that the limiting grid $\Gamma^{(\infty)}$ is optimal or even a local minimum. We refer to the Appendix for a brief discussion and an closed formula for the Hessian of $\mathcal{G}$.
 
\bigskip
The first claim of this  theorem relies on a boundedness assumption for the sequence $\big(\Gamma^{(k)}\big)_{k\ge 0}$. This condition is of course satisfied if  the support of the distribution $\mu$ of $X$  is compactly  supported and  $\mathbf{U}= \overline{{\cal H}({\rm supp}(\mu))}$  (closed convex hull of the support of $\mu$) since we assume that     $\Gamma^{(0)}\subset \mathbf{U}$: then, for every $k\!\in \N$, so will be the case for $\Gamma^{(k)}$ as emphasized in the description of the procedure.

Claim~$(b)$ emphasizes that, by an appropriate choice of (the quadratic quantization error) $\Gamma^{(0)}$ 
may imply the boundedness of the whole sequence of iterates $(\Gamma^{(k)})_{k\ge 0}$. 
This approach known by practitioners as the {\em splitting method}  is investigated  in the next subsection.

\smallskip
We will prove this theorem step by step, establishing intermediary results, often under less stringent assumptions than the above theorem, which may have their own interest.


\subsection{Possibly unbounded inputs: the splitting method}  \label{sec:split}
 In practical computations, if one aims at computing (hopefully) optimal  quantization grids (or $CVT$)  of $X$ on  a wide range of  levels $N$ (see $e.g.$ the website \url{www.quantize.maths-fi.com}),  the so-called {\em splitting method}  appears as an extremely efficient ``level-by-level" procedure. The principle is to compute the grids in a telescopic manner based on their size $N$: assume we have access to an optimal grid of size $N$. Imagine we add to this grid an $(N+1)^{th}$ $\mathbf{U}$-valued component $e.g.$ sampled randomly from the distribution of $X$ (or any $\mathbf{U}$-supported distribution). Doing so, we make up a grid which  has for sure a lower quadratic mean quantization error than any grid with $N$ points. This grid of size $N+1$ is likely to lie in the attracting basin of the optimal  quadratic quantizer (or $CVT$) at level $N+1$. This is often observed in practice and, even if not optimal, it provides at least very good quantizers. Its main the trade-off is that it requires a systematic, hence heavily time consuming, simulation.   Many variants or potential  improvements can be implemented (like adding an optimal quantizer of size $N_0\ge 2$ at each new initialization to directly obtain (hopefully) optimal quantizers of size $N+N_0$ (see $e.g.$~\cite{PAPR}).


\bigskip \noindent\textbf{Splitting Assumption (on the starting grid):}  {\em Let $N\ge 1$ be such that {\rm card}$({\rm supp}\mu)>N$ (where $\mu =\P_{_X}$). Let  $\Gamma_{N-1}^{*}=\{x^*_i, \,i=1,\ldots,n\}\subset \mathbf{U}$   be an optimal grid of size $N-1$ for $\P_{_X}$ and let $x_N^{(0)}\!\in {\rm supp}(\mu) \setminus \Gamma_{N-1}^*$. The Lloyd algorithm is initialized as follows: }
\begin{equation}\label{eq:Gamma0}
\Gamma^{(0)} = \Gamma_{N-1}^{*} \cup \{x_N^{(0)}\} =\Big\{x_1^*, x_2^*, \cdots , x_{N-1}^*, x_N^{(0)}\Big\}.
\end{equation}

One natural way to generate the additional element  $x^{(0)}_{_N}$ is in practice to simulate randomly a copy of $X$ (in fact one can also simulate a copy of a random vector whose distribution is equivalent to that of $X$, see the remark below). 

\bigskip
\noindent {\bf Remark.} If $\mu$ has a density $\rho$, it is more efficient to simulate according to the distribution whose probability density is proportional to $\rho^{\frac{d}{d+2}}$ (one checks that this function is integrable if $X\!\in L^{2+\eta}(\P)$ for an $\eta >0$, see $e.g.$~\cite{Fou}). The reason is that  the resulting distribution $\mu_{d}= \kappa_d \rho^{\frac{d}{d+2}}.\lambda_d $ provides the best random $N$-quantizers of $\mu = \rho.\lambda_d$ for every $N\ge 1$ in the following sense: the asymptotic minimization problem
\[
\inf\Big\{\limsup_{N\to +\infty} N^{\frac 2d}\E \big(\min_{1\le i\le N} |X-Y_i|^2\big),\; Y_1,\ldots,Y_{_N} \; i.i.d., \; \perp\!\!\!\perp X,\; Y_1\sim \nu\Big\}
\]
has $\nu= \mu_{d}$ as a solution ($\perp\!\!\!\perp$ denotes here  independence).

\medskip If the Splitting Assumption on the initialization of the procedure is satisfied,  the grids of the iterations in the Lloyd algorithm share an interesting  property (which implies their global boundedness). The arguments developed in the proof below  are  close  to those used in the proof of the existence of an optimal quantizer at level $N$ (see $e.g.$~\cite{PAG0, Fou}).

\begin{Pro}\label{pro:Npts} Assume  the Splitting Assumption~\eqref{eq:Gamma0}.

\smallskip 
\noindent $(a)$ The quantization error induced by  the  grid $\Gamma^{(0)}$ defined by~\eqref{eq:Gamma0} is strictly smaller than $e_{N-1}(X)$ ($i.e.$  that of the optimal $(N-1)$-quantizer $\Gamma^*_{_{N-1}}$) so that $|\Gamma^{(k)}|=N$ and 
\[
e(\Gamma^{(k)}, X ) = \big \|X - \widehat{X}^{\Gamma^{(k)}} \big\|_2 \le e(\Gamma^{(0)}, X ) < e_{N-1}(X). 
\] 

\noindent $(b)$ The sequence of iterated grids $(\Gamma^{(k)})_{k\ge 0}$ is bounded in $\R^d$ $i.e.$ there exists a compact set $K\subset \R^d$ such that 
\[
\forall\, k\!\in \N,\quad \Gamma^{(k)}\subset K.
\]
\noindent $(c)$  All the limiting values of the sequence of consistent representations   $(x^{(k)})_{k\ge 0}$   have $N$ pairwise distinct components ($i.e.$ the    the sequence of grids $(\Gamma^{(k)})_{k\ge 0}$  has asymptotically full size full size $N$).
\end{Pro}

\noindent
\textbf{Proof}: $(a)$ The (squared) quadratic means quantization error (or ``distortion value") induced by  $\Gamma^{(0)}$ satisfies:
\[
e(\Gamma^{(0)}, X )^2 = \mathcal{G}(x_1^*, x_2^*, \cdots , x_{(N-1)}^*, x_N^{(0)}) = \mathbb{E}\Big( \min_{1\le i \le N} |X - x_i^{(0)}| ^2 \Big).
\]
where $ \, x_j^{(0)} = x_j^*$, $1 \le j \le N-1$. 
Let $\varepsilon = \frac13 \min_{1\le i \le N-1} |x_i^*- x_N^{(0)}|$. For every $\xi \!\in B(x_N^{(0)}, \varepsilon)$, 
\[
|\xi-x_i^{(0)}|\ge | x^*_i-x^{(0)}_{_N}|-|\xi- x^{(0)}_{_N}| \ge 3\varepsilon -\varepsilon= 2\varepsilon > |\xi- x^{(0)}_{_N}|^2.
\]
Consequently, on the event $A_{\varepsilon}= \{X\!\in  B(x_N^{(0)}, \varepsilon)\}$, we have  
\[
\min_{1\le i \le N} |X(\omega) - x_i^{(0)}|_{2}^2 = |X(\omega) - x_N^{(0)}|_{2}^2  >  \min_{1\le i \le N-1} |X(\omega) - x_i^{(0)}|_{2}
\]
whereas we always have that $\min_{1\le i \le N} |X(\omega) - x_i^{(0)}|^2\le \min_{1\le i \le N-1} |X(\omega) - x_i^{(0)}|^2$.  On the  other hand we know that $x^{(0)}_{_N}\!\in {\rm supp}(\P_X)$ so that $\P(A_{\varepsilon})>0$ which implies in turn that 
\[
e(\Gamma^{(0)}, X )^2 < e_{N-1}(X)^2.
\]
One concludes by noting that Lloyd's algorithm makes the sequence of (squared) quadratic quantization error induced by the iterated grids  non-increasing.

\smallskip
\noindent $(b)$ 
Suppose that there exists  a subsequence $(\varphi(k))$ and a component $i_0$ such that  $|x_{i_0}^{(\varphi(k))}|\to +\infty $ as $k\to +\infty$.
Then, by re-extracting finitely many subsequences, we can split the set $\{ 1, \cdots , N\}$ into two disjoint non empty subsets $I$ and $I^c$ and find a subsequence (still denoted $(\varphi(k))_{k\ge 0}$ for convenience) such that 
\[ 
I =  \Big\{ j\!\in \{ 1 \ldots, N\}   \mbox{ such that }\, \big| x_{j}^{(\varphi(k))} \big| \rightarrow +\infty  \Big\} \neq \emptyset
\]
and 
\[ I^c = \Big\{ j\!\in \{ 1 \ldots, N\}   \mbox{ such that }  \,  x_{j}^{(\varphi(k))}   \rightarrow  x_{j}^{\infty} \in \mathbb{R}^d \Big\}.
\]
  It is clear that 
\[
\forall\, \xi \!\in \R^d,\quad \lim_{k \to+ \infty } \min_{1 \le i \le N} \big| \xi - x_{i}^{(\varphi(k))} \big|^2 =  \lim_{k \to \infty } \min_{j \in I^c} \big| \xi - x_{j}^{(\varphi(k))} \big|^2 = \min_{j \in I^c} | \xi - x_{j}^{\infty} |^2.
\]
%
Therefore, it follows from Fatou's Lemma that 
\begin{align}
\liminf_{k \to \infty } \mathcal{G}(\Gamma^{(\varphi(k))}) &= \liminf_{k \to \infty } 
 \mathbb{E} \left( \min_{1 \le i \le N} \big | x - x_{i}^{(\varphi(k))} \big |^2 \right) \notag \\
&\ge  \mathbb{E} \left( \liminf_{k \to \infty }  \min_{1 \le i \le N} \big | x - x_{i}^{(\varphi(k))} \big |^2 \right) \notag \\
&=   \mathbb{E} \left( \liminf_{k \to \infty }  \min_{j \in I^c} \big | x - x_{i}^{(\varphi(k))} \big |^2 \right) \notag \\
&=  \mathbb{E} \left( \min_{j \in I^c} \big | x - x_{j}^{\infty} \big |^2 \right)  \notag \\
&\ge e_{ |I^c|}(X)^2\notag.
\end{align}
Finally $e_{ |I^c|}(X)^2\ge e_{N-1}(X)^2$  since $I$ contains at least $i_0$.  
On the other hand, we know from $(a)$ that $\mathcal{G}(\Gamma^{(\varphi(k))}) \le \mathcal{G}(\Gamma^{(0)}) $, so that we get the following contradictory inequality
\[
e_{N-1}(X)^2\le \liminf_{k \to \infty } \mathcal{G}(\Gamma^{(\varphi(k))}) \le  \mathcal{G}(\Gamma^{(0)}) < e_{N-1}(X)^2. 
\]
Consequently $I$ is empty which completes the proof.

\smallskip
\noindent $(c)$ The proof is similar to that  of the above item $(b)$. Assume $x^{(\varphi(k))}\to x^{(\infty)}$. If there exists   two components  $x_i^{(\varphi(k))}$ and $x_j^{(\varphi(k))}$ converging  to  $x^{\infty}_i= x^{\infty}_j$. Then, owing to Fatou's Lemma 
\[
 e(\Gamma^{(\infty)}, X )^2 =\mathbb{E}\left( \min_{1\le i \le N} |X - x_i^{\infty}|_2^2 \right)  \le \liminf_k  \mathbb{E}\left( \min_{1\le i \le N} |X - x_i^{(\varphi(k))}|_2^2 \right)\le  e(\Gamma^{(\varphi(0))}, X ) < e_{N-1}(X) 
 \]
which yields  a contradiction. $\qquad\Box$

\subsection{Convergence of the Lloyd procedure under a boundedness assumption}
 From now on, our aim is to investigate the structure of the set $\Upsilon^{\infty}$ of limiting grids of the sequence $(\Gamma^{(k)})_{k\ge 0}$ $i.e.$  the set of grids $\Gamma^{(\infty)}$ such that there exists a subsequence $(\varphi(k))_{k\ge 0}$ for which $\displaystyle \lim_{k\to +\infty} \Gamma^{(\varphi(k))}= \Gamma^{(\infty)}$.

The two lemmas and the proposition below establish several properties of the iterated grids $(\Gamma^{(k)})_{k\ge 0}$ which are the basic ``bricks" of the proof of Claim~$(a)$ of Theorem~\ref{thm:main}.

\begin{Lem}\label{lem:nomerge} Let $C$ be a closed convex set of $\R^d$ with non-empty interior and let $\nu$ be a Borel probability distribution such that $\int_{\R^d} |\xi|^2\nu(d\xi) <+\infty$ and $\nu(C)>0$. Furthermore assume that, either $\nu$   satisfies $\nu(\stackrel{\circ}{C})>0$, or $\nu$ assigns no mass to hyperplanes. Then the function defined on $C$ by 
\[
I_C : y\longmapsto \int_C |y-\xi|^2\nu(d\xi) 
\]
is continuous,  strictly convex and atteins its unique minimum  at $y^*_{_C} = \displaystyle \frac{\int_C \xi\nu(d\xi)}{\nu(C)} \in\, \stackrel{\circ}{C}$.
\end{Lem}

\noindent {\bf Proof.}
Elementary computations show that 
\begin{eqnarray*}
I_C (y)&=& \nu(C) |y|^2 +\int_C |\xi|^2\nu(d\xi) -2\left(y|\int_C \xi\nu(d\xi)\right)\\
&=& \nu(C)\left( |y-y^*_{_C}|^2 +\int_C |\xi-y^*_{_C}|^2\frac{\nu(d\xi)}{\nu(C)} \right)
\end{eqnarray*}
so that $ \displaystyle  y^*_{_C}= {\rm argmin}_C I_C$.

\smallskip
Now assume  $y^*_{_C} \!\in \partial C$. There exists a supporting hyperplane $H^*$ to $C$ at $y^*_{_C}$ defined by $H^* = y^*_{_C} +\vec{u}^{\perp}$, $|\vec{u}|=1$. For every $\xi \!\in C$, $(\xi-y^*_{_C}|\vec{u})\ge 0$ so that 
\[
0= (0|\vec{u})= \frac{1}{\nu(C)}
 \int_C \underbrace{(\xi-y^*_{_C}|\vec{u})}_{\ge 0}\nu(d\xi)
 \]
so that $(\xi-y^*_{_C}|\vec{u})=0$ $\nu(d\xi)$-$a.s.$ $i.e.$ $\xi\!\in H^*$ $\nu(d\xi)$-$a.s.$. This leads to a contradiction with the assumptions made on $\nu$.\hfill $\quad\Box$

\begin{Lem} \label{lem:cvgce} Assume that a subsequence $\Gamma^{(\varphi(k))}\to \Gamma^{(\infty)}$ as $k\to +\infty$ and that the boundary of the Voronoi tessellation of $\Gamma^{(\infty)}$ are $\P$-negligible. Let  $Y \in L^1(\Omega,{\cal A}, \mathbb{P})$. Then 
\[ 
\mathbb{E} ( Y | \widehat{X}^{\Gamma^{(\varphi(k))}} ) \rightarrow \mathbb{E} ( Y | \widehat{X}^{\Gamma^{(\infty)}} ) \qquad a. \, s. \; \mbox{ as } \; k\to +\infty.
\]
\end{Lem}

\noindent \textbf{Proof of Lemma~\ref{lem:cvgce}}: By definition,
\begin{align}
 \mathbb{E} ( Y | \widehat{X}^{\Gamma^{(\varphi(k))}} ) 
 &= \sum_{i=1}^N \frac{\E\, Y  \mathbf{1}_{X \in C_i(\Gamma^{(\varphi(k))})}  } {\mathbb{P}(X \in C_i(\Gamma^{(\varphi(k))}))}\mbox{\bf 1}_{\{X \in C_i(\Gamma^{(k)})\}}. \notag
\end{align}

 We know from what precedes that, since  $\Gamma^{(\varphi(k))} \rightarrow \Gamma^{(\infty)}$, if $X(\omega)\!\in \mathbb{R}^d \setminus \bigcup_{i=_1}^N \partial C_i(\Gamma^{(\infty)})$, the functions
\[
 \mathbf{1}_{X \in C_i(\Gamma^{(\varphi(k))})} (\omega) Y (\omega)  \rightarrow \mathbf{1}_{X \in C_i(\Gamma^{(\infty)})} (\omega) Y (\omega)\quad\mbox{ as } \; k\to +\infty,\; i=1,\ldots,N,
 \]
Our assumption on $\Gamma^{(\infty)}$ implies that these convergences hold $\P$-$a.s.$.
We conclude by Lebesgue's dominated convergence theorem that  for every $i\!\in \{1,\ldots,N\}$ 
\[ 
\E \big(Y(\xi)^a \mathbf{1}_{X \in C_i(\Gamma^{(\varphi(k))})}   \big)\longrightarrow \E \big(Y(\xi)^a \mathbf{1}_{X \in C_i(\Gamma^{(\infty)})}\big)
\]
Therefore, with our convention for the index $i$ such that $\mathbb{P}(X \in C_i(\Gamma^{(\infty)}))=0$ (if any), we get by applying the above convergence to $Y$ and ${\bf 1}$
\begin{align}
\mathbb{E} ( Y | \widehat{X}^{\Gamma^{(\varphi(k))}} ) &= \sum_{i=1}^N \frac{\E \big(Y \mathbf{1}_{X \in C_i(\Gamma^{(\varphi(k))})} \big)}{\mathbb{P}(X \in C_i(\Gamma^{(\varphi(k))}))} \mbox{\bf 1}_{\{X \in C_i(\Gamma^{(k)})\}}\notag \\
 &\xrightarrow[k\to\infty] \, \sum_{i=1}^N \frac{\E \big(Y \mathbf{1}_{X \in C_i(\Gamma^{(\infty)})}\big)}{\mathbb{P}(X \in C_i(\Gamma^{(\infty)})}\mbox{\bf 1}_{\{X \in C_i(\Gamma^{(\infty)})\}} =  \mathbb{E} ( Y | \widehat{X}^{\Gamma^{(\infty)}} ) \qquad \P\mbox{-}a. s. \qquad\qquad\hfill \Box \notag
\end{align}

\begin{Pro}[Grid convergence~I] \label{prop:cvgceI} Assume  that $\Gamma^{(0)}$ is a $\mathbf{U}$-valued grid and that the iterates $(\Gamma^{(k)})_{k\ge 0}$ of  the Lloyd algorithm are bounded.


\smallskip
\noindent $(a)$  Assume that $\mu$ assigns no mass to hyperplanes  and ${\rm supp}(\mu) = \mathbf{U}$ ($i.e.$ is convex). 
If the sequence $(\Gamma^{(k)})_{k\ge 0}$ of iterations of the Lloyd procedure is bounded ($e.g.$ because $\mathbf{U}$ is itself bounded) then
\[
\liminf_k \min_{i\neq j} |x^{(k)}_i-x^{(k)}_j|>0
\]
$i.e.$ no components of the grids get asymptotically stuck as $k$ goes to infinity.

%

%
%

\smallskip
\noindent $(b)$ Assume $X\!\in L^2(\P)$. Let $\Gamma^{(\infty)}$ be a limiting grid of $(\Gamma^{(k)})_{k\ge 0}$. If  the boundary of  the Voronoi cells of $\Gamma^{(\infty)}$ are $\P_{_X}$-negligible $i.e.$  $\P\big( X \!\in \cup_i \partial C_i(\Gamma^{(\infty)})\big) =0$, then the  grid $\Gamma^{(\infty)}$ is stationary $i.e.$ it is a fixed point of the Lloyd map $T_{_N}$ or equivalently that (any of)  its induced Voronoi tessellation(s) is a CVT. Moreover if $\Gamma^{(\varphi(k))}\to \Gamma^{(\infty)}$, then 
\[
\widehat X^{\Gamma^{(\varphi(k))}}\stackrel{a.s.\,\&\, L^2}{\longrightarrow} \widehat X^{\infty}\quad\mbox{ as } \; k\to +\infty.
\]

\noindent $(c)$ If the distribution of $X$ assigns no mass to hyperplanes, then $\nabla \mathcal{G}(\Gamma^{(k)})\to 0$ as $k\to + \infty$.

\smallskip
\noindent $(d)$ If the distribution of $X$ 
has a convex support $\mathbf{U}={\rm supp}(\mu)$, then the consistent representations of $(\Gamma^{(k)})_{k\ge 0}$ satisfy
\[
\sum_{k\ge 0}  |x^{(k+1)} - x^{(k)} |^2_{(\R^d)^N} <+\infty. 
\]
In particular, $x^{(k+1)}-x^{(k)}\to 0$ as $k\to +\infty$. Hence, the set $\mathcal{X}_{\infty}$ of (consistent representations of) limiting grids  of $(\Gamma^{(k)})_{k\ge 0}$ is a (compact) connected subset of $\mathbf{U}^N$.
%
\end{Pro}

%
\noindent {\bf Remarks.} $\bullet$ In the literature the convergence of the gradient  of the iterated grids as established  in $(c)$   is sometimes known as ``weak convergence"  of the Lloyd procedure.

\smallskip 
\noindent $\bullet$ Lemma~\ref{lem:nomerge}  improves a result obtained in~\cite{Maria} which show that component do not asymptotically merge in Lloyd's procedure even when the Splitting Assumption is not satisfied. In~\cite{Maria}, the distribution $\mu$ is supposed to be absolutely continuous  with a density a continuous $\rho$, everywhere  strictly positive on its support.  Our method of proof allows for relaxing this absolute continuity assumption. 
%
%


\bigskip 
\noindent\textbf{Proof of Proposition~\ref{prop:cvgceI}}: $(a)$ 
For every $i,\, j\!\in \{1,\ldots,N\}$, $i\neq j$, we define the median hyperplane of $x^{(k)}_i$ and $x^{(k)}_j$ by 
\[
\vec{u}^{(k)}_{ij}= \frac{x^{(k)}_i-x^{(k)}_j}{|x^{(k)}_i-x^{(k)}_j|},\;  H^{(k)}_{ij} = \frac{x^{(k)}_i+x^{(k)}_j}{2} +\vec{u}_{ij}^{\perp} \quad \mbox{ and }\quad \]
We define together the affine form    
\[
\varphi_{ij}(\xi) = \Big(\xi- \frac{x^{(k)}_i+x^{(k)}_j}{2} \Big|\vec{u}_{ij}  \Big),\quad \xi\!\in \R^d
\]
which satisfies $\varphi^{(k)}_{ij} \ge 0$ on $C_i(\Gamma^{(k)})$ and $\varphi^{(k)}_{ij}(x^{(k)}_i) =\frac 12  |x^{(k)}_i-x^{(k)}_j|>0$. Moreover, note that $\varphi^{(k)}_{ji}=-\varphi^{(k)}_{ij}$ and that
\[
\bigcap_{j\neq i} \Big\{\varphi_{ij}^{\infty}> 0\Big\}= \stackrel{\circ}{C}_i(\Gamma^{(k)})\subset C_i(\Gamma^{(k)})\subset \overline{C}_i(\Gamma^{(k)})\subset \bigcap_{j\neq i} \{\varphi_{ij}^{\infty}\ge 0\}.
\] 
The sequence of iterated grids $(\Gamma^{(k)})_{k\ge 0}$ being bounded by assumption, we may assume without loss of generality, up to an extraction $(\varphi(k))_{k\ge 0}$, 
\[
x^{(\varphi(k))}_i\to x^{\infty}_i,\; i\!\in \{1,\ldots,N\},\quad \vec{u}^{(\varphi(k))}_{ij}\to u^{\infty}_{ij}, \; i,\,j\!\in \{1,\ldots,N\}
\]
and
\[
x^{(\varphi(k)+1)}_i \to \tilde x^{\infty}_i, \; i\!\in \{1,\ldots,N\}\quad\mbox{ as } k\to+\infty.
\]
Set for every $i,\, j\!\in \{1,\ldots,N\}$
\[
\varphi_{ij}^{\infty} = \lim_{\varphi(k)\to \infty}\varphi^{(\varphi(k))}_{ij} \quad \mbox{and}\quad C^{\infty}_i = \bigcap_{j\neq i} \{\varphi_{ij}^{\infty}\ge 0\}.
\]
Hence,  for every $i\!\in \{1,\ldots,N\}$, $C^{\infty}_i$ is a closed polyhedral convex set containing $x_i$. It also contains $\tilde x_i$ since $x^{(\varphi(k)+1)}_i\!\in C_i(\Gamma^{(\varphi(k))})$ owing to the stationarity property. Then, for every $j\!\in \{1,\ldots,N\}$, 
\[
\varphi_{ij}^{\infty }(\tilde x^{\infty}_i) =\lim_{\varphi(k)\to \infty} \varphi_{ij}^{(\varphi(k))}(x^{(\varphi(k)+1)}_i)\ge 0
\]
since $\varphi^{(\varphi(k))}_{ij}$  uniformly on compact sets toward $\varphi^{\infty}_i$ (simple convergence of affine forms in finite dimension implies locally uniform convergence).

Assume there exists $i_0\!\in \{1,\ldots,N\}$ such that  $I_0=\{i\,|\,  x^{\infty}_i= x^{\infty}_{i_0}\}$  is not reduces to $\{i_0\}$   $i.e.$ contains at least two indices. 

It is clear that  
\[
\emptyset \neq \Big\{\xi\!\in \R^d\,|\, |\xi-x^{\infty}_{i_0}|<d\big(\xi, \Gamma^{(\infty)}\setminus\{x_{i_0}\}\big)\Big\}\subset \bigcup_{i\in I_0} C^{\infty}_i.
\]
The above nonempty open set has non zero $\mu$-mass since $x_{i_0}\!\in \mathbf{U}$. Hence, there exists two indices $i_1,i_2\!\in I_0$ such that $\mu(C^{\infty}_{i_1})$ and $\mu(C^{\infty}_{i_2})>0$. First note that  $\varphi^{\infty}_{i_1i_2}(x^{\infty}_{i_1})\ge 0$ and $\varphi_{i_2i_1}(x^{\infty}_{i_2})\ge 0$ but both quantities being opposite since $x^{\infty}_{i_1}=x^{\infty}_{i_2}$ they are equal to $0$ which means that 
\[
x^{\infty}_{i_0}=x^{\infty}_{i_1}=x^{\infty}_{i_2} \!\in \partial C^{\infty}_{i_1}\cap \partial C^{\infty}_{i_2}.
\]

Since these sets are polyhedral and $\mu$ assigns no mass to hyperplanes, $\mu(\stackrel{\circ}{C}_{i_1})$ and $\mu(\stackrel{\circ}{C}_{i_2})>0$. Then one checks that 
\[
\mbox{\bf 1}_{\stackrel{\hskip -0.25 cm \circ}{C_{i_\ell}^{\infty}} }= \lim_k \mbox{\bf 1}_{\stackrel{\hskip -0.75 cm _\circ}{C_{i_\ell}^{(\varphi(k))}} }, \; \ell=1,2
\]
so that, still using that $\mu$ assigns no mass to the boundaries of these polyhedral convex sets,  we get 
\[
0<\mu( \stackrel{\hskip -0.25 cm _\circ}{C_{i_{\ell}}^{\infty}})= \lim_{k\to+\infty}\mu\big(\stackrel{\hskip -0.75cm _\circ}{C_{i_\ell}^{(\varphi(k))}}\big), \; \ell=1,2.
\]
Set $\varepsilon_0= \min_{\ell=1,2} \mu(\stackrel{\hskip -0.25 cm\circ}{C_{i_\ell}^{\infty}})>0$. It follows form~(\ref{eq:ldeltagrid}) that
\[
\| X-\widehat X^{\Gamma{(\varphi(k))}}  \|^2_2-\|   X-\widehat X^{\Gamma{(\varphi(k)+1)}}  \|^2_2\ge \min_{\ell=1,2} \mu\big(\stackrel{\!\!\!\!\circ}{C_{i_\ell}}^{(\varphi(k))}\big) \big(\big|x^{(\varphi(k)+1)}_{i_1}-x^{(\varphi(k))}_{i_1}\big|^2+ \big|x^{(\varphi(k)+1)}_{i_2}-x^{(\varphi(k))}_{i_2}\big|^2\big).
\]
Letting $\varphi(k)$ go to infinity implies that 
\[
\varepsilon_0\big(\big|\tilde x^{\infty}_{i_1}-x^{\infty}_{i_1}\big|^2+ \big|\tilde  x^{\infty}_{i_2}-x^{\infty}_{i_2}\big|^2\big)\le 0
\]
since $\| X-\widehat X^{\Gamma{(k)}}\|_2$ is a converging sequence. Hence $\tilde x^{\infty}_{i_{\ell}}=x^{\infty}_{i_{\ell}}$, $\ell=1,2$, which in turn implies that $\tilde x^{\infty}_{i_1}= \tilde x^{\infty}_{i_1}=x_{i_0}$.

One shows likewise, still taking advantage of the $\mu$-negligibility of the boundaries of the polyhedral sets $C^{\infty}_i$,  that
\[
\lim_{k\to+\infty} \int_{C_{i_{\ell}}(\Gamma^{(\varphi(k))})}\xi \mu(d\xi)=  \int_{C_{i_{\ell}}^{\infty}}\xi \mu(d\xi),\; \ell=1,2.
\]
Passing to the limit in the stationary equation satisfies by $x_{_{\ell}}^{(\varphi(k))}$, $\ell=1,2$,    finally implies that 
\[
\tilde x^{\infty}_{i_{\ell}} = \frac{ \int_{C_{i_{\ell}}^{\infty}}\xi \mu(d\xi)}{\mu(C_{i_{\ell}}^{\infty})}\!\in\, \stackrel{\hskip -0.25 cm \circ}{C_{i_{\ell}}^{\infty}},\quad\ell=1,2.
\]
But, owing to Lemma~\ref{lem:nomerge} (see also~\cite{Fou}, p.22),  this implies that $\tilde x^{\infty}_{i_{\ell}} \!\in\, \stackrel{\hskip -0.25 cm \circ}{C_{i_{\ell}}^{\infty}}$, $\ell=1,2$. This yields a contradiction to the fact that both $x_{\ell}$ are equal $x_{i_0}$.

\smallskip
\noindent  $(b)$ Let  $\xi \! \in \mathbb{R}^d \setminus \bigcup_{i=_1}^N \partial C_i(\Gamma^{(\infty)})$. 
Then $\xi$  belong to the interior of one of the tessels $C_i(\Gamma^{(\infty)})$, say  $\mathring{C}_{i_0}( \Gamma^{(\infty)} )$. Hence,    $|\xi- x_{i_0}^{\infty}|$ is strictly smaller than $\min_{i\neq i_0}|\xi-x_i^{\infty}|$. Consequently, there exists   an $n(\xi)\in\mathbb{N}^*$ such that for all $n\ge n(\xi)$, 
\[ 
|\xi - x_{i_0}^{(\varphi(k))}| < \min_{i \neq i_0 } |\xi  - x_{i}^{(\varphi(k))}|
\]
or equivalently $\xi \! \in \mathring{C}_{i_0}( \Gamma^{(\varphi(k))} )$. 
Thus,  for every $\xi\! \in \mathbb{R}^d \setminus \bigcup_{i=_1}^N \partial C_i(\Gamma^{(\infty)})$,  
\[
\operatorname{Proj}_{\Gamma^{(\varphi(k))}} (\xi) =  \sum_{i=1}^N x_i^{(\varphi(k))} \mathbf{1}_{C_i(\Gamma^{(\varphi(k))})}(\xi) \xrightarrow[k\to \infty] \,  \sum_{i=1}^N x_i^{\infty} \mathbf{1}_{C_i(\Gamma^{(\infty)})}(\xi ) = \operatorname{Proj}_{\Gamma^{(\infty)}} (\xi ).
\]
This clearly implies that $\P(d\omega)$-$a.s.$, $\widehat{X}^{\Gamma^{(\varphi(k))}}(\omega) \to  \widehat{X}^{\Gamma^{(\infty)}}(\omega)$ as $k\to+ \infty$ 
%
%
%
%
since $\mathbb{P}( X\! \in \bigcup_{i=_1}^N \partial C_i(\Gamma^{(\infty)}))=0$. To carry on the proof, we need the following lemma. 
%
%
%
%
%

%

\noindent If we set $Y=X$ in Lemma~\ref{lem:cvgce}, then  
\[ 
\mathbb{E} ( X | \widehat{X}^{\Gamma^{(\varphi(k))}} ) \longrightarrow \mathbb{E} ( X | \widehat{X}^{\Gamma^{(\infty)}} )  \qquad a.s. 
 \]
Moreover,  the  sequence $\big(\E(X\,|\, \widehat X^{\Gamma^{(k)}})\big)_{k\ge 0}$ being $L^2$-uniformly integrable  since $X\!\in L^2(\P)$, the above convergence also holds in $L^2(\P)$.
 
Now, since $ \varphi(k+1) \ge \varphi(k)+1$ and $\widetilde X^{\Gamma^{(\varphi(k)+1)}}= \E \big(X\,|\, \widehat X^{\Gamma^{(\varphi(k))}}\big)$ is $\Gamma^{(\varphi(k)+1)}$-valued, we derive from Lemma~\ref{Lem:seminal} that 
\[ 
\forall k \in \mathbb{N}^*, \quad \| X -  \widehat{X}^{\Gamma^{(\varphi(k+1))}} \|_2\le  \|  X -  \widehat{X}^{\Gamma^{(\varphi(k)+1)}} \|_2 \le\|X-   \widetilde X^{\Gamma^{(\varphi(k)+1)}} \|_{_2}=   \|  X -   \E \big(X\,|\, \widehat X^{\Gamma^{(\varphi(k))}}\big) \|_2.
 \]
Since we know that  $\widehat{X}^{\Gamma^{(\varphi(k))}} \longrightarrow  \widehat{X}^{\Gamma^{(\infty)}}$ $\P$-$a.s.$, it follows from  from Fatou's Lemma that
\[
 \|  X -  \widehat{X}^{\Gamma^{(\infty)}} \|_2  \le   \liminf_{k \to \infty } \| X -  \widehat{X}^{\Gamma^{(\varphi(k+1))}} \|_2.
\]
On the other hand,  we derive from the convergence $\mathbb{E} ( X | \widehat{X}^{\Gamma^{(\varphi(k))}} ) \longrightarrow \mathbb{E} ( X | \widehat{X}^{\Gamma^{(\infty)}} )$ in $L^2(\P)$ that  
\[
 \lim_{k \to \infty }\|  X -   \E \big(X\,|\, \widehat X^{\Gamma^{(\varphi(k))}}\big) \|_2 =   \|  X -  \mathbb{E} ( X | \widehat{X}^{\Gamma^{(\infty)}} )  \|_2.  
\]
%
%
so that
\[
 \|  X -  \widehat{X}^{\Gamma^{(\infty)}} \|_2  \le  \|  X -  \mathbb{E} ( X | \widehat{X}^{\Gamma^{(\infty)}} )  \|_2
\]
which in turn  implies  by the very definition of conditional expectation as an orthogonal projection on $L^2(\widehat X^{\Gamma^{(\infty)}})$ that  
\[  
\widehat{X}^{\Gamma^{(\infty)}} =  \mathbb{E} ( X | \widehat{X}^{\Gamma^{(\infty)}} ) \qquad \P\mbox{-}a.s.
\]

\smallskip 
\noindent $(c)$  First we note that, for any grid $\Gamma=\{x_1,\ldots,x_{_N}\}$, Schwarz's Inequality implies 
\begin{align}
|\nabla \mathcal{G}(\Gamma)|_2^2 &= \sum_{i=1}^N \left\vert \frac{\partial \mathcal{G}}{\partial x_i}\left(\Gamma \right)\right\vert^2 
&
= 4 \sum_{i=1}^N \left\vert \int_{C_i(\Gamma} \left(  x_i  - \xi \right) \mathbb{P}(\mathrm{d} \xi) \right\vert^2  
\le  4 \sum_{i=1}^N \int_{C_i(\Gamma}\left\vert x_i - \xi \right\vert^2 \mathbb{P}(\mathrm{d} \xi)  
= 4 \mathcal{G}(\Gamma) \notag
\end{align}
so that the sequence $\big(\nabla \mathcal{G}(\Gamma^{(k)})\big)_{k\ge 0}$ is bounded. Now, as  $\P_X$   assigns no mass to the boundary of any Voronoi tessellations,   we derive from what precedes that any limiting grid of $(\Gamma^{(k)})_{k\ge 0}$ is stationary $i.e.$ $\nabla \mathcal{G}(\Gamma^{(\infty)})=0$. It is clear by an extraction procedure that $0$ is the only limiting value for the bounded sequence $\big(\nabla \mathcal{G}(\Gamma^{(k)})\big)_{k\ge 0}$ which consequently converges toward $0$.

\smallskip 
\noindent $(d)$   The set of limiting grids of the sequence $(\Gamma^{(k)})_{k\ge 0}$ is compact by construction. Its connectedness will classically follow  from 
\[ 
|\Gamma^{(k+1)} - \Gamma^{(k)} |  \xrightarrow[]{k\to \infty} 0
\] 
 (see~$e.g.$~\cite{LP}). 
 
 It is clear from  item~$(a)$ that $K_{\Gamma}=\{ \Gamma^{(k)},\, k\ge 0\}\cup \Upsilon_{\infty}$ is a compact whose intersection with the closed set $\{(x_1,\ldots,x_{_N})\!\in (\R^d)^N,\, x_i\neq x_j,\, i\neq j\}^c$. Hence, $K_{\Gamma}$ stands at a positive distance of this set  or, equivalently,  there exists $\delta>0$, such that 
\[
 \forall\, k\ge 0,\quad \forall\, i,\, j\in \{1,\ldots,N\}, \, i\neq j, \quad |x^{(k)}_i-x^{(k)}_j|\ge \delta.
\]
As a consequence,  every Voronoi cell satisfies
\[
B^{\hskip -0.18 cm ^{^\circ}}(x^{(k)}_i, \delta/2) \subset C_i(\Gamma^{(k)})
\]
where $B^{\hskip -0.18 cm ^{^\circ}} (\xi, r)$ denotes the  open ball with center $\xi$ and radius $r$. Since $\mathbf{U}={\rm supp}(\mu)$ is a (closed) convex set and $\Gamma^{(0)}\subset {\rm supp}(\mu)$, then $\Gamma{(k)}\subset \mathbf{U}$ for every $k\!\in \N$. As a consequence,  $K_{\Gamma}\subset \mathbf{U}$ which in turn implies that the function $\xi\mapsto \P\big(X\!\in B^{\hskip -0.18 cm ^{^\circ}} (\xi, \delta/2)  \big)$ is (strictly) positive on the compact set $K_{\Gamma}$, so that we can define
\begin{equation}\label{eq:mstar}
m_* := \inf_{\xi \in \mathbf{K_{\Gamma}}}\P\big(X\!\in B^{\hskip -0.18 cm ^{^\circ}} (\xi, \delta/2)  \big)>0
\end{equation}
which implies in particular that 
\[ 
M^{(k)}_i =  \mathbb{P}(X\!\in C_i(\Gamma^{(k)})) \geqslant m_*,\; 1 \le i \le N,  \; k\ge 0.
\]

Now we introduce the  {\em energy gap}
 $$
 \Delta(k)= \big\|X-\widehat X^{\Gamma^{(k)}} \big\|^2_2-\big\|X-\widetilde X^{\Gamma^{(k+1)}} \big\|^2_2
 $$ 
 known to be non-negative by~\eqref{eq:step1} in Lemma~\ref{Lem:seminal}. Then
\begin{eqnarray}
\Delta(k) &=& \sum_{j=1}^N \int_{C_j(\Gamma^{(k)})} | \xi - x_j^{(k)}|^2 \mathbb{P}_{_X}(\mathrm{d} \xi) - \sum_{j=1}^N \int_{C_j(\Gamma^{(k)})} | \xi - x_j^{(k+1)}|^2 \mathbb{P}_{_X}(\mathrm{d} \xi). \notag  \\ 
&=&  \sum_{j=1}^N \int_{C_j(\Gamma^{(k)})} (| x_j^{(k)} |^2 - | x_j^{(k+1)}|^2) + 2\left (  x_j^{(k+1)} - x_j^{(k)} | \xi\right ) \mathbb{P}_{_X}(\mathrm{d} \xi)  \notag  \\ 
&=&  \sum_{j=1}^N M^{(k)}_j(| x_j^{(k)} |^2 -  |x_j^{(k+1)}|^2)  + 2 \Big( x_j^{(k+1)} - x_j^{(k)} | \int_{C_j(\Gamma^{(k)})}  \xi \,\mathbb{P}_{_X}(\mathrm{d} \xi) \Big) \notag  \\
&=&  \sum_{j=1}^N M^{(k)}_j(| x_j^{(k)} |^2 -  | x_j^{(k+1)}|^2)  + 2  M^{(k)}_j\left ( x_j^{(k+1)} - x_j^{(k)} |  x_j^{(k+1)}  \right)\notag  \\ 
&=& \sum_{j=1}^N M^{(k)}_j| x_j^{(k)} -  x_j^{(k+1)}|^2 \notag \\
&\ge & m_* |\Gamma^{(k+1)} - \Gamma^{(k)} |^2_{(\R^d)^N}.\label{eq:ldeltagrid}
\end{eqnarray}

On the other hand  $\Delta(k)\le \big \|X-\widehat X^{\Gamma^{(k)}}\big\|^2_2-\big\|X-\widehat X^{\Gamma^{(k+1)}} \big\|^2_2$ so  that, finally,
\[
\qquad \qquad \qquad \sum_{k\ge 0}  |\Gamma^{(k+1)} - \Gamma^{(k)} |^2_{(\R^d)^N}\le \frac{1}{m_*} \big\|X-\widehat X^{\Gamma^{(0)}}\big\|^2_2<+\infty.\qquad\qquad  \qquad \qquad\qquad \qquad  \hfill \Box
\]

\medskip
%
%

%
\noindent {\sc Comments.} $\bullet$ The restrictions on the possible limiting grids in  Theorem~\ref{thm:main}   do not imply uniqueness in general in higher dimensions: the symmetry properties shared by the distribution itself already induces multiple limiting grids as emphasized by the case of the  multivariate normal distribution $\mathcal{N}(0, I_{d})$.
In fact, for an orthogonal  matrix $P\! \in \mathcal{O}_{d}$ that is $PP^* = I_{d}$ ($P^*$ stands for the transpose of $P$) and any optimal  grid $\Gamma=\{x_1,\ldots,x_{_N}\}$,
\[
\|X-\widehat X^{P\Gamma}\|^2_2 = \E \min_{1\le i\le N} |X-Px_i|^2 = \E \min_{1\le i\le N} |P^*X-x_i|^2= \E \min_{1\le i\le N}|X-x_i|^2= \|X-\widehat X^{\Gamma}\|^2_2.
\]

\noindent $\bullet$ For  distribution with less symmetries  like $X \sim \mathcal{N}(0, \Sigma_{d})$ with $\Sigma = \operatorname{diag} (\lambda_1, \lambda_2 , \cdots, \lambda_d), \, 0< \lambda_1 < \lambda_2 , \cdots < \lambda_d$, one can reasonably hope that at least local uniqueness of stationary grids (CVT) holds true. One way to check that is to establish that the Hessian  $D^2\mathcal{G}(\Gamma)$ is invertible at each stationary grid $\Gamma$. A closed form is available for this Hessian  (see~$e.g.$~\cite{FOPA}).

 \bigskip
 \noindent
 \textbf{Proof of Theorem~\ref{thm:main}}:  The preliminary claim follows by induction from the structural stationary properties of the iterates and Lemma~\ref{lem:nomerge}.
 
 \smallskip
  \noindent  $(a)$ Combining the results obtained in the above proposition and the convergence $\big\|X-\widehat X^{\Gamma^{(k)}}\big\|_2$ toward a non-negative real number $\ell_{\infty}$ completes the proof. 
 
 \smallskip
 \noindent $(b)$ follows from Proposition~\ref{pro:Npts} in Section~\ref{sec:split} which implies that the sequence $(\Gamma^{(k)})_{k\ge 0}$ is bounded. Then one concludes by $(a)$. $\quad \Box$

 \section{A bounded variant of Lloyd's procedure based on spatial estimation of the optimal quantizers}\label{sec:3}
So far  our results are based on the hypothesis that we initialize  the Lloyd algorithm using a grid of size $N$ whose induced quadratic quantization error is lower than the minimal quantization error achievable with a grid of size at most $N-1$. This is clearly the key point  to ensure that the iterates of the procedure remain  bounded. From a practical point of view, this choice for  the initial grid  is not very realistic, in particular if we are processing a ``splitting method": nothing ensures, even if the Lloyd procedure converges at a level $N-1$, that the limiting grid will be optimal with the consequence that  the initialization at level $N$ ``below" $e_{N-1}(X)$ becomes impossible. 

However, we know from theoretical results on optimal vector quantization where the optimal quantizers are located {\em a priori}.   

\subsection{A priori bounds for optimal quantizers}

The following proposition can be found in~\cite{Fou}.

\begin{Pro}   Let $L_{N,X}(c) = \{ \Gamma \,:\, |\Gamma|= N \mbox{ and } e(\Gamma,X) \leqslant c \}$.
Let $c \!\in (0,e_{N-1}(X)]$. There exists $R\!\in (0,+\infty)$ such that 
$$
L_{N,X}(c)\subset B(m_{_X},R) \mbox{ with } m_{_X} = \E\,X.
$$

An upper bound $S$ satisfies the following conditions:
\begin{itemize}
\item[(i)] $\exists\, r>0$ such that $\mathbb{P}\big(X\!\in B(m_{_X} ,r)\big)>0$ and $(\frac{R}{5}-r)^2 \mathbb{P}\big(X\!\in B(m_{_X} ,r)\big)> c$, 
\item[(ii)] $4\displaystyle \int_{B(0,\frac{2R}{5})^{c}} |\xi-m_{_X} |_2^2 \mathbb{P}_{X}(\mathrm{d} \xi) < e_{N-1}(X) - c$.
\end{itemize} 
\end{Pro}
If we specify $c = e_{N }(X)$, then the set  $L_{N,X}(c)$ will be the set of grids corresponding to optimal $N$-quantizers, and  $R$ will be an upper bound of the optimal $N$-quantizers. 

\bigskip Consequently, in order to given a numerical estimation of $R$, we have to estimate the asymptotic behaviour of $e_{N-1}(X) - e_{N }(X)$. We know from~(\cite{Sagna}, see also~\cite{GRLUPA3}) that
 $$
 e_{N-1}(X) - e_N(X) \propto N^{-\frac{d+2}{d}}
 $$
  If $X$ has a standard Gaussian distribution ${\cal N}(m; I_d)$, then 
\[ 
\int_{B(m_{_X} ,\frac{2R}{5})^{c}} |\xi-m_{_X}|^2 \mathbb{P}_X(\mathrm{d} \xi)\int_{B(0 ,\frac{2R}{5})^{c}} |\xi|^2 \mathbb{P}_{_X}(\mathrm{d} \xi)  \propto R^d e^{-\frac{R^2}{2}}.\]
Thus
\[ 
R
 \propto \sqrt{\ln(N)}.
\]  

\medskip More precise results  can be found  in~\cite{Sagna} and~\cite{JULU} for various families of distributions with exponential or polynomial tails at infinity. Under certain condition on $X$, the asymptotic behaviour of $R=R(N)$ can be analyzed sharply as $N$ goes to infinity. Typically if $X\sim {\cal N}(m; I_d)$, 
\[
\lim_N \frac{R(N) }{\sqrt{\log(N) }}= \frac{1}{\sqrt{2}}\Big(1+\frac 2d\Big)^{\frac 12}.
\]  
\subsection{A variant of Lloyd's algorithm}

  Since we know that all  optimal $N$-quantizers are  constrained in a bounded domain (depending on $N$ in  a more or less controlled way), a natural idea is to constrain the exploration of the Lloyd iterates inside it to take advantage of this information. 

\medskip We can fix an area that we are sure that the optimal quantizer have its grids in it. Once an iteration of the algorithm includes some points that go beyond the area, we will be sure that it is not the optimal grid. So we can do something to pull the points back into the area while keeping the error non-increasing.

To this end, we now compute the difference made on the second phase of the Lloyd iteration in condition that we pick a point other than the mass center point. If we take $x_j^{\prime}$ instead of $ x_j^{(k+1)}$ for a certain $j$, the resulting difference in the $j$\textsuperscript{th} Voronoi region created by the second phase will be:

\begin{align} \Delta_j^{\prime}(k) :=& \;\E |X-x^{(k+1)}_j|^2\mbox{\bf 1}_{\{X\in C_j(\Gamma^{(k)}) \}}-   \E |X-x'_j|^2\mbox{\bf 1}_{\{X\in C_j(\Gamma^{(k)}) \}}\notag\\
=&\;\int_{C_j(\Gamma^{(k)})} | \xi - x_j^{(k)}|^2 \mathbb{P} (\mathrm{d}\xi) -\int_{C_j(\Gamma^{(k)})} | \xi - x_j^{\prime}|^2 \mathbb{P} (\mathrm{d}\xi)\notag \\ 
=&\; \int_{C_j(\Gamma^{(k)})} | \xi - x_j^{(k)}|^2 \mathbb{P} (\mathrm{d}\xi) -  \int_{C_j(\Gamma^{(k)})} | \xi - x_j^{(k+1)}|^2 \mathbb{P} (\mathrm{d}\xi) \notag \\ 
 &+ \int_{C_j(\Gamma^{(k)})} | \xi - x_j^{(k+1)}|^2 \mathbb{P} (\mathrm{d}\xi) -\int_{C_j(\Gamma^{(k)})} | \xi - x_j^{\prime}|^2 \mathbb{P} (\mathrm{d}\xi) \notag \\ 
=& \;M_j(k)\Big(| x_j^{(k)} -  x_j^{(k+1)}|^2-| x_j^{\prime} -  x_j^{(k+1)}|^2\Big). \notag 
\end{align}

This shows that {\em if we pick a point which lies at  the same distance to $x^{(k+1)}_j$ as  $x^{(k)}_j$}, then the above difference  becomes zero. The idea is then to choose the point inside the prescribed domain if $x^{(k+1)}_j$ is outside. It is always possible, $e.g.$ by keeping $x^{(k)}_j$ still (although this is probably not the optimal way to proceed).

With this idea we present a modified version of Lloyd's algorithm, which  does not need   the Splitting Assumption, to be run successfully. We set a $R>0$ that all optimal quantizers are in the ball $B(0, R)$. The algorithm is as follows (assuming that $X$ is centered for convenience):

\bigskip
\noindent {\bf Algorithm 2 (Modified Lloyd's algorithm)}:

\smallskip
\noindent {\em Inputs}:
\begin{itemize}
\item  $B(0, R)$, the domain of interest (with $R$ close to $R(N)$) hopefully. 
\item  $\mu=\P_{_X}$ a simulatable probability distribution (with a convex support $\mathbf{U}$) and assigning no mass to hyperplanes.
\item  $\Gamma^{(0)} = \{x_i^{(0)},\, i=1,\ldots,N \}$, the initial set of $N$ generators (starting grid).
\end{itemize}

\noindent 
{\em Pseudo-script}:

\smallskip
\noindent  $\rhd$ At the $k$\textsuperscript{th} iteration:
\begin{enumerate}
\item Compute the position of $\Gamma^{(k+1)} =  \{x_i^{(k)+1},\,i=1,\ldots,N \}$, the mass centroid of $\{C_i^{(k)},\,i=1,\ldots,N \}$. If there is an index $j$ such that  $x_j^{(k+1)} $ lies outside  $B(0, R)$ (and for every such point),  replace the current value  $x_j^{(k+1)} $  by a point  $x_j^{(k+1)\prime}$ defined $e.g.$ by  $x_j^{(k+1)\prime}=\partial B(x_j^{(k+1)} , | x_j^{(k)} -  x_j^{(k+1)}|)\cap [x^{(k)}_j,x^{(k+1)}_j]$ (other choices are possible like choosing it randomly on $\partial B(x_j^{(k+1)} , | x_j^{(k)} -  x_j^{(k+1)}|)\cap B(0,R)$. 
\item If every point lies in the ball $B(0, R)$, then take these mass centroids of $\{C_i^{(k)},\,i=1,\ldots,N \}$ as the new set of generators $\Gamma^{(k+1)} =  \{x_i^{(k+1)},\,i=1,\ldots,N \}$,
\item Construct the Voronoi tessellation $\mathcal{C}(\Gamma^{(k+1)}) =\{C_i^{(k+1)},\,i=1,\ldots,N  \}$ of $\mathbf{U}$ with the grid of generators $\{x_i^{(k+)},\,i=1,\ldots,N \}$.
\end{enumerate}

\noindent  $\rhd$  Repeat the above iteration until a stopping criterion is met. And the output is the CVT $\{C_i^{(n)},\,i=1,\ldots,N \}$ with generators $\{x_i^{(n)} ,\,i=1,\ldots,N \}$ in $\mathbf{U}$.

\medskip
\noindent  $\rhd$  end.

\begin{figure}[h!]
   \centering
    \includegraphics[width=1.08\textwidth]{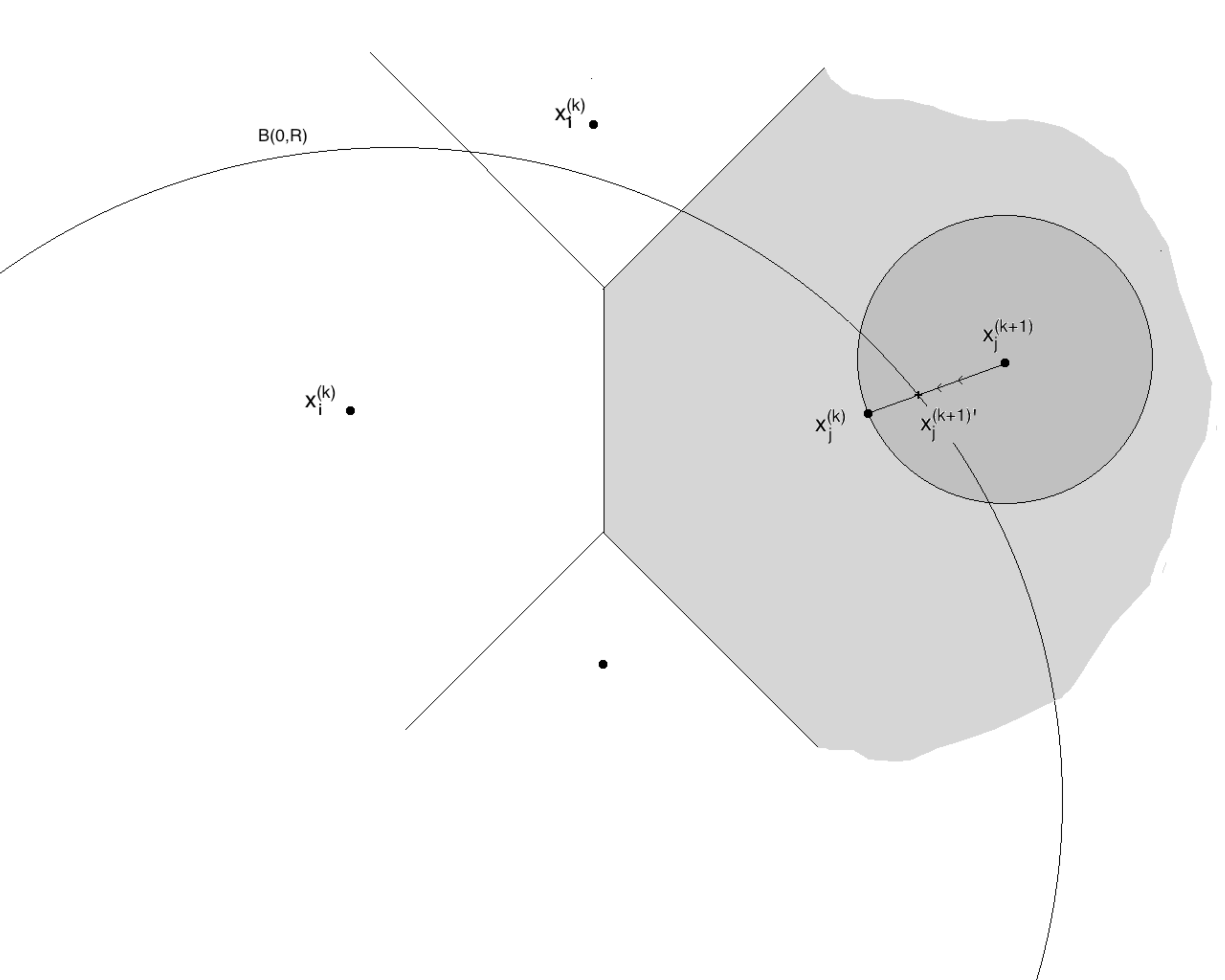}
   \caption{Modified Lloyd's procedure:   $x^{(k+1)}_j$ is replaced by $x^{(k+1)'}_j= [x^{(k)}_j, x^{(k+1)}_j ]\cap\partial B(0,R)$}
  \label{fig:lloyd_modif}
  \end{figure}

\medskip
In this new version  we  modify Phase~I~(\ref{eq:step1}) of the   Lloyd procedure in such a way that the quadratic approximation error $\|X-\widetilde X^{\Gamma^{(k+1)}}\|_2$ still decreases.
The second phase being unchanged, so this modified Lloyd  algorithm is still   energy descending   and furthermore it lives in the ball $B(0, R)$. The trade-off is that with the $R$ fixed at the beginning of the algorithm, we loose the (theoretical) possibility that the iterated sequence cruises very far during the iterations to  finally come  back  with a lower energy. Therefore the energy level of the new limit points will be higher than setting $R=+\infty$. We also see that the larger the radius $R$ we take, the lower  limit energy level we can get. 

\medskip Another trade-off of the modified procedure is that it does not guarantee Lemma 3 because we do not use the Splitting Assumption. In this case we cannot prove the non-degeneracy of the limit grid by this global energy reasoning. However, we can now rely on Proposition~\ref{lem:nomerge} to ensure that no merging occurs.


 \medskip
 \noindent  {\bf Provisional remarks.}  One verifies on numerical implementation of Lloyd algorithm, that the main default that slows down the procedure is more the freezing of one component of the grid which is`` too far from the core" of the support of the distribution $\mu$ than the explosion of the grid with components going to infinity. but in some sense these  seeming radically different behavior are the two sides of the same coin and the above procedure is an efficient way to prevent these parasitic effects. Though, in practice we proceed in a less formal way: using the theoretical estimates on the radius of the distribution allows for an adequate choice of the initial grid $\Gamma^{(0)}$ as confirmed by various numerical experiments carried $e.g.$ in~\cite{SagnaPhD}.

\section{Appendix: Numerical detection of the nature of a limiting grid} 
We provide here a formula for the Hessian of the distortion value  function $\mathcal{G}$ when the distribution $\mu$ of $X$ is absolutely continuous with density function $\rho$. From such a formula it is possible, at least numerically in low dimensions, to detect the status of a stationary grid/quantizer in terms of stability: local minimum, saddle point, etc. 

As a first step we need a Lemma which provides a formula for the differentiation of integrals overs teasels of the Voronoi partition of an $N$-tuple $x=(x_1,\ldots,x_{_N})$. 
\begin{Lem} Let $\varphi\!\in {\cal C}(\R^d, \R)$. Set for every $x\!\in \R^d$ with pairwise distinct components, 
\[
\Phi_i(x) = \int_{C_i(x)}\varphi(\xi)\lambda_d(d\xi), \; i=1,\ldots,N,
\]
where $\lambda_d$ denotes the Lebesque measure on $(\R^d, {\cal B}or(\R^d))$. Then $\Phi_i$ is continuously differentiable on the open set of $N$-tuples with  pairwise distinct components and 
\[
\forall\, j\!\in \{1,\ldots,N\},\; j\neq i,\; \frac{\partial \Phi_i}{\partial x_j}(x)= \int_{\overline{C}_i(x)\cap \overline{C}_j(x)}\varphi(\xi) \left( \frac 12 n_{x}^{ij} + \frac{1}{|x_i-x_j|} \left(\frac{x_i+x_j}{2}-\xi\right)\right)\lambda_x^{ij}(d\xi)
\]
where  $\lambda_x^{ij}(d\xi)$ denotes the Lebesque measure on the median hyperplane  $H_{ij}^x$ of $x_i$ and $x_j$ and $n^{ij}_x = \frac{x_j-x_i}{|x_j-x_i|}$. Furthermore,
\[
\frac{\partial \Phi_i}{\partial x_i} (x) = -\sum_{j\neq i}\frac{\partial \Phi_j}{\partial x_i}(x).
\]
\end{Lem}
We refer to~\cite{FOPA} or~\cite{POL} for  a proof. This leads to the announced general result concerning the Hessian of the distortion function $\mathcal{G}$. We set for every $u= (u^1,\ldots,u^d),\,v=(v^1,\ldots,v^d)\!\in \R^d$, $u\otimes v= [u^iv^j]_{1\le i,j\le d}$ and $I_d=[\delta_{ij}]_{1\le i,j\le d}$ ($\delta_{ij}$ denotes the Kronecker symbol).


\begin{Pro} Let  $\mu =\P_{_X} =  \rho.\lambda_d$ with $\rho$ continuous. Then,  for every $i,\, j\!\in \{1,\ldots,N\}$, $i\neq j$, 
\[
\frac{\partial^2 \mathcal{G}}{\partial x_i\partial x_j}(x)=  
\int_{\overline{C}_i(x)\cap \overline{C}_j(x)}(x_i-\xi)\otimes \left(\frac 12 n^{ij}_x +\frac{1}{|x_i-x_j|}\Big(\frac{x_i+x_j}{2}-\xi\Big)\right)\rho(\xi)\lambda^{ij}_x(d\xi) \quad\mbox{if }i\neq j
\]
and 
\[
\frac{\partial^2 \mathcal{G}^{\ell}}{\partial x_i\partial x_j}(x)= \mu(C_i(x)) I_d+ \sum_{j\neq i} \int_{\overline{C}_i(x)\cap \overline{C}_j(x)}(x_j-\xi)\otimes \left(\frac 12 n^{ij}_x -\frac{1}{|x_i-x_j|}\Big(\frac{x_i+x_j}{2}-\xi\Big)\right)\rho(\xi) \lambda^{ij}_x(d\xi)  \; \mbox{if }\,i = j.
\]

\end{Pro}

This formula is used $e.g.$ in~\cite{FOPA} to show the instability of ``square", "hyper"-rectangular stationary grids for the uniform distribution over the unit hypercube for Kohonen's Self-Organizing Maps ($SOM$). When the neighborhood function of the $SOM$ is degenerated (no true neighbor) the $SOM$ amounts to the $CLVQ$ and its equilibrium points are those of the Lloyd procedure, with the same (un-)stability properties. These quantities also appear in the asymptotic variance of  the $CLT$ established in~\cite{POL}.
\end{document}